\documentclass[12pt,thmsa]{article}
\usepackage{amsfonts}
\usepackage{amssymb}
\newtheorem{teo}{Theorem}[section]

\newtheorem{lema}[teo]{Lemma}

\newtheorem{prop}[teo]{Proposition}

\newtheorem{obs2}[teo]{Remark}

\newtheorem{tea}{Theorem}[subsection]

\newtheorem{no2}[teo]{Note}

\newtheorem{no3}[tea]{Note}

\newcommand{\Gal}{{\rm Gal}}
\newcommand{\Frob}{{\rm Frob }}

\newcommand{\F}{{\mathbb{F} }}
\newcommand{\Q}{{\mathbb{Q} }}
\newcommand{\mod}{{\rm mod}}

\newcommand{\PSL}{{\rm PSL}}
\newcommand{\SL}{{\rm SL}}
\newcommand{\PGL}{{\rm PGL}}

\newcommand{\GL}{{\rm GL}}

\newcommand{\Image}{{\rm Image}}

\newcommand{\PGSp}{{\rm PGSp}}

\newcommand{\GSp}{{\rm GSp}}

\title{Uniform behavior of families of Galois
representations on Siegel modular forms and the Endoscopy Conjecture}
\author{Luis V. Dieulefait \\ Departament d'Algebra i Geometria, Universitat de Barcelona
 \thanks{
  Research supported by project MTM2006-04895, MECD, Spain}} 

\begin{document}

\maketitle
\begin{abstract}
 We
prove the following uniformity principle: if one of the Galois representations in
the family attached to a genus two Siegel cusp form of weight $k>3$, ``semistable" and
with multiplicity one,
is reducible (for an odd prime $p$),
then all the representations in the family are reducible. This, combined with Serre's conjecture (which is now a theorem) gives
 a proof of the Endoscopy Conjecture.

\end{abstract}

\section{Introduction}
In this article, we will consider a genus two Siegel modular form $f$
of level $N$ and weight $k>3$ (and multiplicity one)
and the family of four dimensional symplectic Galois representations
attached to it. We assume also that we are in a case where this family is ``semistable". In [D1], we have treated the level $1$ case, giving conditions on $f$ to ensure that
these Galois representations have generically large image. In
particular we have imposed an irreducibility condition on one
characteristic polynomial of Frobenius (see [D1], condition (4.8)) to obtain a large image
result. Furthermore, with the same irreducibility condition, we showed
in [D2], again for the level $1$ case, that for every $p > 4k-5$ the $p$-adic representations are
irreducible. The only possible reducible case to be considered is
the case of two $2$-dimensional irreducible components having the same
determinant (all other cases can not occur if $f$ is not of Saito-Kurokawa type,
cf. [D1], [D2]), and from the results of [D2] this case can only
happen 
if all characteristic polynomials are reducible, i.e., the
$2$-dimensional components will have coefficients in the same field that the
$4$-dimensional representations: the field $E$
generated by the eigenvalues of $f$. \\

In section 2 we will
generalize the main results of [D1] and [D2] to the semistable case.\\

One of the consequences of Tate's conjecture on the Siegel threefold  is
that reducibility for the Galois representations
attached to $f$ must be a uniform property: if it is verified at one
prime, then all the representations in the family are reducible.
 In this article, we will prove this uniformity principle:

\begin{teo}
\label{teo:unifo} : Let $f$ be a genus $2$ Siegel cuspidal Hecke eigenform
of weight $k >3$ and level $N$,
 having
multiplicity one, such that the attached Galois representations $\rho_{f,\lambda}$ are ``semistable". Suppose that for some odd prime $\ell_0 \nmid N$, $\lambda_0 \mid
\ell_0$, the representation $\rho_{f,\lambda_0}$ is reducible. Then the
representations $\rho_{f,\lambda}$ are reducible for every
$\lambda$.\\
Moreover, if this happens, either $f$ is of Saito-Kurokawa type or $f$ is endoscopic.
\end{teo}

After excluding the Saito-Kurokawa case, we will prove the  result more generally for
compatible families of geometric, pure and symplectic four-dimensional Galois
representations which are ``semistable".\\

A previous version of this preprint dates from 2003, and since several papers using the results contained there have appeared since then, we prefer to present first (sections 2 and 3) the results contained in that early version, which constitute the ``core" of this paper: this corresponds to the proof of ``uniformity of reducibility" for ``almost every" prime, and with the extra condition $\ell_0 > 4k-5$. At the end of the paper (section 4) we will indicate how to (easily) remove this assumption on $\ell_0$. Finally, we will prove that a standard combination of these results with Serre's conjecture allows us to remove the ``almost every" in the result and gives also the Endoscopy Conjecture, i.e., the modularity (up to twist) of the irreducible components.\\ 

What follows is a brief description of the tools that will appear in the proofs in the ``core" part.
We will use (as in [D2], section 4) as starting point Taylor's recent results on
the Fontaine-Mazur conjecture and the meromorphic continuation of $L$-functions
 for odd two-dimensional Galois
representations (see [T2], [T3] and [T4]). Then, we will combine some of the results and
 techniques in [D1] (in particular the information about the description of
 the action of inertia obtained via $p$-adic Hodge theory) with Ribet's
results (see [R])
on two-dimensional semistable Galois representations (slightly
generalized to higher weights), and finally Cebotarev density theorem,
the fundamental theorem of Galois theory, and some
group theory will suffice for the proof.\\

\section{Preliminaries}
 
  As we already explained, the goal of sections 2 and 3 is to prove a theorem which is weaker than theorem \ref{teo:unifo}, namely we will prove the following:
  
  \begin{teo}
  \label{teo:unifow} Let $f$ be a genus $2$ Siegel cuspidal Hecke eigenform
of weight $k >3$ and level $N$,
 having
multiplicity one, such that the attached Galois representations $\rho_{f,\lambda}$ are ``semistable". Suppose that for some prime $\ell_0 > 4k-5, \ell_0 \nmid N$, $\lambda_0 \mid
\ell_0$, the representation $\rho_{f,\lambda_0}$ is reducible. Then the
representations $\rho_{f,\lambda}$ are reducible for almost every
$\lambda$.\\
  
  \end{teo}

From now on we will make the following assumption: $f$ is a genus $2$
level $N$ Siegel cuspidal Hecke eigenform of weight $k>3$,
 having multiplicity one, and
not of Saito-Kurokawa type (theorem \ref{teo:unifow} is trivial in the Saito-Kurokawa
 case, where by construction the Galois representations are reducible, with one $2$-dimensional and two $1$-dimensional components). Let $E = \mathbb{Q}( \{ a_n \})$  be the field generated by its Hecke eigenvalues. Then,
 there is a compatible family of Galois representations constructed by
 Taylor [T1] and Weissauer [W2] verifying the following:\\
  For any prime
number $\ell$ and any extension $\lambda$ of $\ell$ to $E$ we have
a continuous Galois representation
$$ \rho_{f,\lambda} : G_{\mathbb{Q}} \rightarrow \GSp(4, \overline{E}_\lambda)$$
unramified outside $\ell N$ and with characteristic polynomial of
$ \rho_{f,\lambda}(\Frob \; p)) $ equal to
$$ Pol_p (x) = x^4 - a_p x^3 + (a_p^2 - a_{p^2} - p^{2k-4})x^2 - a_p p^{2k-3} x
+ p^{4k-6}$$
for every $p \nmid \ell N$. If $ \rho_{f, \lambda} $ is absolutely
irreducible, then it is defined over $E_\lambda$.\\

In general, we can not guarantee that the field of definition is
$E_\lambda$, but the residual representation $\bar{\rho}_{f,\lambda}$
can be formally defined in any case (see [D1])
as a representation  defined over the residue field of $\lambda$,
$\mathbb{F}_\lambda$. Nevertheless, not knowing the field of definition
of the representations that we will study is not a serious problem, we can
work instead with the ``field of coefficients" (*), i.e., the field
generated by the coefficients of the characteristic polynomials
$Pol_p(x)$, this field contains all the information we need.\\

The representations $\rho_{f,\lambda}$
 are known to have the following properties (cf [W1], [W2],
[D1]): they are pure (Ramanujan conjecture is satisfied) and if $\ell \nmid N$
they are crystalline with Hodge-Tate weights $\{ 0, k-2, k-1, 2k -3\}$. This
 last property makes possible, via Fontaine-Laffaille theory, to obtain
 a precise description of the action of the inertia group at $\ell$ on
 the residual representation $\overline{\rho}_{f,\lambda}$: it acts
 through fundamental characters of level one or two, with exponents
 equal to the Hodge-Tate weights (see [D1] for more details). \\

We will need a further restriction: we want the representations to be ``semistable" at every prime of $N$ (see the definition below).\\

Since we will not use the fact that our representations are modular, we can change to the more general setting of a family of
 four-dimensional symplectic Galois representations $\{\rho_\lambda \}$
  with coefficients in a number field
 $E$ (not necessarily defined over $E_\lambda$, see (*)),
 $\det{\rho_\lambda} = \chi^{4k-6}$, which are pure,
  and such that there exists a finite set $S$ with,
 for every $\ell \not\in S$, $\rho_\lambda$ unramified outside $\{\ell\} \cup
 S$, crystalline at $\ell$ with Hodge-Tate weights as above, and
 ``semistable" at primes in $S$, i.e., verifying the following:
 $\rho_\lambda$ restricted to $I_q$ is a unipotent group for every $q \in
 S$.\\
For every $p \not\in S$ we still denote $Pol_p (x)$ the characteristic
polynomial of the image of $\Frob \; p$ and $a_p$ the trace of this
image. The representations being symplectic, we have the standard
factorization
  $$  Pol_p(x) = (x^2- (a_p/2 + \sqrt{d_p}) x + p^{2k-3})
  (x^2- (a_p/2 - \sqrt{d_p}) x + p^{2k-3}) \qquad \; (2.1) $$
  The results of generically large image and irreducibility
   proved in previous articles for the level $1$ case
  (see [D1], theorem 4.2, and [D2], theorems 2.1 and 4.1)
 hold also in this generality:

 \begin{teo}
 \label{teo:large}: Let $\{ \rho_\lambda \}$ be a family of Galois
 representations verifying the above properties, with $k>3$.
Assume that there is a prime $p \not\in S$
 such that
$$ \sqrt{d_p} \notin
  E \qquad \qquad (2.2) $$
  where $d_p$ is defined by formula (2.1).
 Then for all but finitely many of the primes verifying
$$  d_p \notin
 (\mathbb{F_\lambda})^2  $$
 and, more generally, for
  all primes $\lambda$ in $E$ except at most for a set of Dirichlet
density $0$, the
 image of $\rho_{\lambda}$ is
 $$ A_{\lambda}^k = \{ g \in \GSp( 4, \mathcal{O}_{E_\lambda}) :
\det(g) \in ( \mathbb{Q}_\ell^*)^{4k-6} \} ,$$
 where $\mathcal{O}_{E_\lambda}$ denotes the
 ring of integers of $E_\lambda$.\\
Keeping condition (2.2) we also have: for every
prime $\ell \geq 4k-5, \ell \not\in S$, $\lambda \mid
\ell$, the representation $\rho_\lambda$ is absolutely
irreducible.
\end{teo}

Differences with the level 1 case:\\
 The proof of the above results given in [D1] and [D2] extends automatically
 to
 the semistable case:
recall that the determination of the images is done by
considering the image of the residual $\mod \:
\lambda$  representations and eliminating all
non-maximal proper subgroups of $\GSp(4,
\mathbb{F}_\lambda)$. When considering
reducible cases (cf. [D1], sections 4.1 and 4.2) if we
allow  arbitrary ramification at a finite set $S$ then
we have to allow  the character appearing as
one-dimensional component or determinant of a
two-dimensional component of a reducible
$\bar{\rho}_\lambda$ to ramify at $S$, but in the
semistable case it is easy to see that this character will not ramify
at primes in $S$.
The same applies to the case of image
equal to a group $G$ having a reducible index $2$ normal
subgroup $M$ (cf. [D1], section 4.4), the quadratic
Galois character $G/M$ can not ramify at primes of $S$ if we
assume semistability. Up to these easy remarks, all the
proof translates word by word to the semistable
case. \\

Remark 1: Recall that condition (2.2) was introduced
(cf. [D1]) specifically  to deal with the case where the image of
$\bar{\rho}_\lambda$ is reducible, with two
$2$-dimensional irreducible components of the same
determinant. All other cases of non-maximal image can be discarded,
for almost every prime, without using condition (2.2).
\\

Remark 2: In [D1], the large images result was
proved (for the case of conductor 1) with an
additional condition, called ``untwisted"
: this condition was imposed to eliminate the
possibility that the projective residual image falls
in a smaller symplectic group $\PGSp(4, k')$, $k'$ a proper
subfield of $k$, where
$k$ is the field generated by the
traces of the residual representation. We have not
included a similar condition in the above theorem because
 in the following lemma, we will
explain that this condition is superfluous, i.e., that
the case of smaller projective symplectic group can never happen
if we assume semistability. In particular, this
applies to level 1 Siegel cusp forms, so  the
condition ``untwisted" can be removed from theorem 4.2
of [D1].

\begin{lema}
\label{teo:untwisted}: Let $\{ \rho_\lambda \}$ be a
compatible families of Galois representations as
above (in particular, a semistable family). Then for
every prime $q > 2k-2, q \not\in S$ and $Q$ a prime in
$E$ dividing $q$,
if we call $G$ the image of $\bar{\rho}_Q$ and $P(G)$
its projectivization, $P(G)$ lies in $\PGSp(4, k)$ if
and only if $G$ lies in $\GSp(4, k)$, for every subfield
$k$ of $\mathbb{F}_Q$.

\end{lema}

Proof: A similar result, for semistable
two-dimensional representations, is lemma 2.4 in [R].
The proof given there translates word by word, once we
have explained why in our case we also have an element
$c$ in the inertia group $I_q$ such that $\chi(c)$ is
a generator of $\mathbb{F}_q^*$ and the trace of
$\bar{\rho}_Q (c)$ is a non-zero element of
$\mathbb{F}_q$ ( we know a priori, from the
description of the action of $I_q$, that this trace
will be in $\mathbb{F}_q$, what requires a proof is
the fact that it is not $0$).\\
We have given in [D1], proposition 3.1, a description
of the action of $I_q$ that applies in the current
situation, because we are assuming that $\rho_Q$ is symplectic and
crystalline with Hodge-Tate
 weights $\{ 0, k-2, k-1, 2k-3 \}$, and $q>2k-2$. Let $\psi$ be a level 2 fundamental
character, and take $c \in I_q$ such that $\psi(c)$
generates $\mathbb{F}_{q^2}^*$.  We have four
possibilities for the trace of $c$, whose values are,
after a suitable factorization:
$$ (1 + \chi(c)^{k-1} ) ( 1 +  \chi(c)^{k-2}) $$
$$  (\psi(c)^{k-2} + \psi(c)^{ (k-2) q} )
        ( \psi(c)^{k-1}  + \psi(c)^{ (k-1) q })   $$
$$ ( 1 + \psi(c)^{  (k-2)+(k-1) q } )
           (1 + \psi(c)^{ (k-1)+(k-2) q } ) $$
$$ (  \psi(c)^{k-2} + \psi(c)^{(k-1) q}  )
           (\psi(c)^{k-1} + \psi(c)^{(k-2) q}  ) $$
In all cases, the inequality $q> 2k-2$ implies that
these traces are not $0$.\\

\section{ Uniformity of reducibility}
At this point, we can say that the validity or not
of condition (2.2) at some prime $p \not\in S$
determines the behavior of the family of
representations $\rho_\lambda$: If condition (2.2) is
satisfied, then we have generically large image and
irreducibility for every $\ell$ sufficiently large
compared with the weights. \\ What happens if
condition (2.2) is not satisfied at any prime? This implies that
the factorization (2.1) takes place over $E$, i.e.,
that for every $p \not\in S$, $Pol_p(x)$ reduces over
$E$.
The coefficients of all characteristic polynomials
$Pol_p(x)$ generate an order $\mathcal{O}$ of $E$, and
if we restrict to primes $\lambda$ not dividing the
conductor of this order (we are neglecting only
finitely many primes), we see that the field generated
by the coefficients of the  $\mod \; \lambda$ reduction
of all the $Pol_p(x)$ gives the whole
$\mathbb{F}_\lambda$. Thus, we see that for almost
every prime, the failure of (2.2) implies that
$\bar{\rho}_\lambda$ has its image in $\GSp(4,
\mathbb{F}_\lambda)$ and not in a smaller symplectic,
but all characteristic
polynomials reduce over $\mathbb{F}_\lambda$:
in this case the image can not be the whole symplectic
group, because in the group $\GSp(4, \mathbb{F}_\lambda)$ most of the matrices have IRREDUCIBLE characteristic polynomial, and we know that for almost every prime only
one possibility  (see  remark 1 after theorem \ref{teo:large},
 and lemma \ref{teo:untwisted}) remains:

\begin{lema}
\label{teo:residualreducible}: Let $\{ \rho_\lambda \}
$ be as in the previous section, and assume that for
every $p \not\in S$, condition (2.2) is not satisfied.
Then, for almost every prime $\lambda$, the residual
representation
$\bar{\rho}_\lambda$ is  reducible with two
$2$-dimensional irreducible components of the same
determinant.
\end{lema}

\subsection{A reducible member in the family: Residual consequences}
From now on,
 assume that for a prime $q > 4k-5 , q \not\in S
$, $Q \mid q$, the $Q$-adic representation $\rho_Q$ is
reducible, we know (using semistability and purity)
that the only possible case is the case of two
$2$-dimensional irreducible components both with
determinant $\chi^{2k-3}$. Thus
  we have:
$$ \rho_Q \cong \sigma_{1,Q} \oplus \sigma_{2,Q} \quad \quad (3.1)$$
 Since this representations is reducible, the last part of 
theorem \ref{teo:large} implies that condition (2.2) must fail at every
prime. Therefore, $\sigma_{1,Q}$ and $\sigma_{2,Q}$
will also have coefficients in $E$ and lemma \ref{teo:residualreducible}
implies that for
 every prime $\lambda$ in a cofinite set $\Lambda$ of primes of $E$,
  $\bar{\rho}_\lambda$ will verify:
$$\bar{\rho}_\lambda \cong \pi_{1,\lambda} \oplus \pi_{2, \lambda}$$
where $\pi_{i,\lambda}$ is an irreducible two dimensional
representation defined over $\mathbb{F}_\lambda$ having
determinant $\chi^{2k-3}$, for $i=1,2$ and for every $\lambda \in
\Lambda$.\\
Moreover, we can determine the image of $\bar{\rho}_\lambda$ for almost
every prime in $\Lambda$:
\begin{lema}
\label{teo:reducibleperogrande}: Keep the  above assumptions. For
every prime $\lambda \in \Lambda_2$, a cofinite subset of $\Lambda$,
 the image of $\bar{\rho}_\lambda$ is
a  subgroup of $\GSp(4, \mathbb{F}_\lambda)$ conjugated to
$M_\lambda = \{A \times B \in \GL(2,\F_{1,\lambda}) \times \GL(2, \F_{2,\lambda}) :
\det(A) = \det(B) \in \F_\ell^{2k-3} \}$, where
 $\F_{1,\lambda}, \F_{2,\lambda} \subseteq \F_\lambda$ are the fields of coefficients
 of $\pi_{1,\lambda} $ and $\pi_{2, \lambda}$ .
\end{lema}

Proof: We have assumed that the representations $\rho_\lambda$ have a
finite ramification set $S$ and they are semistable at every prime $q \in
S$. A fortiori, the same applies to their residual components
$\pi_{i,\lambda}$. Moreover, these two dimensional representations are
irreducible for every $\lambda \in \Lambda$. In a similar situation,
Ribet has proved a large image result for $\ell \geq 5$, but he assumes
that the action of $I_\ell$, given by fundamental characters of level
$1$ or $2$, has weights (i.e., exponents of the fundamental characters)
$0$ and $1$. The main point of his proof is to exclude the dihedral case.
In our case, using the information on the Hodge-Tate
decomposition, we have this extra condition at $I_\ell$ also verified by the twisted
representation $\pi_{i,\lambda} \otimes \chi^{-k+2} $ for, say, $i=2$
(cf. [D1],[D2]). On the other hand, for $\pi_{1,\lambda}$
 Ribet's result still holds if we
restrict to primes $\ell > 4k-5$, because the weights of the action of $I_\ell$
 being $0$ and $2k-3$, the projectivization of the image if $I_\ell$
 gives a cyclic group of order $(\ell \pm 1) / \gcd(\ell \pm 1 , 2k-3) >
 2$, and this is all that you need to follow Ribet's argument. We also have a
  statement as lemma \ref{teo:untwisted} for these two dimensional representations,
   again adapting lemma 2.4 in [R].\\
We conclude (cf. [R], theorem 2.5 and the remark after)
 that for $\ell$ sufficiently large, the images of both irreducible
 components are conjugated to the subgroup of matrices in $\GL(2, \F_{i,\lambda})$
 with determinant in $\F_\ell^{2k-3}$.\\
 Finally, to prove that the image of $\bar{\rho}_\lambda$ is as we
 want, it remains to show that the Galois fields corresponding to
 $P(\pi_{1,\lambda})$ and $P(\pi_{2,\lambda})$ are disjoint
  ($P$ denotes projectivization). These fields having Galois groups isomorphic to
   the simple groups $\PGL(2, \F_{i,\lambda})$ or $\PSL(2, \F_{i,\lambda})$,
    they are either
   disjoint or equal: the second is not possible because the
   restriction of these two projective representations to $I_\ell$ are different,
   and this proves the result.\\

\subsection{A reducible member in the family: $\lambda$-adic consequences}
In the decomposition (3.1) of $\rho_Q$ it is clear that $\sigma_{1,Q}$
has Hodge-Tate weights $\{0, 2k-3 \}$ and $\sigma_{2,Q}$ has Hodge-Tate
weights $\{ k-2, k-1 \}$ (or viceversa).\\ Now, we invoke a result of
Taylor (see [T2] and [T3], recall that $q > 4k-5$)
 asserting that for a representation such as $\sigma_{1,Q}$ it is possible to
  find a totally
  real number field $F$ such that it is modular when restricted to this
  field,
   and therefore it
   agrees on $F$ with the
   $Q$-adic motivic irreducible
   Galois representation (constructed by Blasius and Rogawski)
    attached to
   a Hilbert modular form $h$.
    This implies that $\sigma_{1, Q}$ appears in the
 cohomology of the restriction of scalars of the motive $M_h$ associated to
 $h$, and it can be checked from the fact that the $Q$-adic
 representation of the absolute Galois group of $F$
 attached to $h$ has descended to a $2$-dimensional
 representation of $G_\Q$, Cebotarev density theorem, and the fact that
 all modular Galois representations in the family $\{ \sigma_{h,\lambda} \}$
  attached to $h$ are
 known to be irreducible, that the whole family descends to a compatible family
 $\{ \sigma_{1,\lambda}  \}$ of
 Galois representations of $G_\Q$
 containing $\sigma_{1,Q}$. To do this, one has to write the representation $\sigma_{1,Q}$
 as in the proof of theorem 6.6 in [T3], and define the representations $\sigma_{1,\lambda}$
  formally in the same way using the strongly compatible families associated to the
  base change of $h$ to each $E_i$ (recall that, for each $i$, $F/E_i$ is soluble, cf.
   [T3]). Then, following an idea suggested to us by R. Taylor, one can check that the
   virtual representations $\sigma_{1,\lambda}$ constructed this way are true Galois
    representations by applying the arguments of [T4], section 533.\\
It follows from the main result of [T3] that  the family
     $\{  \sigma_{1,\lambda }\}$
      is a strongly compatible family (cf. [T3] for the definition) of
  Galois representations. Strong compatibility
   proves the last steps of the following:
  \begin{prop}
  \label{teo:Taylor} Let $\rho_Q$ be as above, reducible as in (3.1),
  and let $\sigma_{1,Q}$ be its irreducible component having Hodge-Tate
  weights $\{0, 2k-3 \}$. Then, there exists a  compatible family of
  Galois representations $\{ \sigma_{1,\lambda} \}$ containing $
  \sigma_{1,Q}$, such that
   for every $\ell \not\in
  S$, $\lambda \mid \ell$, the representation $\sigma_{1,\lambda}$ is
  unramified outside $\{ \ell \} \cup S$, is crystalline at $\ell$ with
  Hodge-Tate weights $\{ 0 , 2k-3 \}$, and is semistable at every prime
  of $S$. Of course, these representations are pure because $\rho_Q$
  is.
  \end{prop}
  Recall that, the representation $\rho_\lambda$ being
  symplectic, for every $g \in G_\mathbb{Q}$ the roots of $\rho_\lambda (g)$
   come in reciprocal pairs: $\{ \alpha , \chi^{2k-3}(g)/ \alpha ,
    \beta , \chi^{2k-3}(g)/ \beta \}$.\\
  The following lemma is a first approach to compare the
  representations $\sigma_{1,\lambda}$ and $\rho_\lambda$:
\begin{lema}
\label{teo:mismasraices} For every $\ell \not\in S$, $\lambda \mid
\ell$, and every $g \in G_\mathbb{Q}$, the roots of $\sigma_{1,\lambda}(g)$
 form a pair of reciprocal roots of those of $\rho_\lambda (g)$.
\end{lema}
Proof: From the compatibility of the families $\{ \sigma_{1,\lambda } \}$
and $\{  \rho_\lambda \}$
and the fact that $\sigma_{1,Q}$ is a component of $\rho_Q$ the lemma
is obvious for the dense set of Frobenius elements at unramified
places. Then, by continuity and Cebotarev the lemma follows for every
element of $G_\mathbb{Q}$.\\

Recall that $\Lambda_2$ denotes the cofinite set of primes of $E$ where
lemma \ref{teo:reducibleperogrande} is satisfied. We will shrink again
this set by eliminating a finite set of primes, namely, those primes
where the image of $\sigma_{1,\lambda}$ fails to be maximal: in fact,
if we call $E' \subseteq E$ the field of coefficients of this family
of representations and $\mathcal{O}'$ its ring of integers,
 using semistability and again the slight
modification of the methods of [R] to higher weights (as we did before to obtain
lemma \ref{teo:reducibleperogrande}) we see that for almost every prime
$\lambda \in E$ the image of ${\sigma}_{1,\lambda}$ can be
conjugated to the
subgroup of $\GL(2, \mathcal{O}'_\lambda)$ of matrices with determinant
in $\mathbb{Z}_\ell^{2k-3}$ (after proving the similar result
 for the residual representations, we apply a lemma of Serre in [S1]
 that shows that the $\lambda$-adic image is also large).\\
 Remark: Here we need to know that the residual representations
  $\bar{\sigma}_{1,\lambda}$ are
 almost all of them irreducible. This follows again from the good properties
  of the $\lambda$-adic family: purity, the fact
  that they are all crystalline with Hodge-Tate
 weights $\{ 0, 2k-3\}$ (and the uniform description of inertia that one gets from
 this), and semistability.\\

 Thus, we exclude from $\Lambda_2$ the finite set of primes where the
 image of $\sigma_{1,\lambda}$ fails to be maximal, and we obtain a
 cofinite set $\Lambda_3$ where the residual image of $\rho_\lambda$ is
 the full $M_\lambda$ and the image of $\sigma_{1,\lambda}$ is maximal.

We want to extract more information from the relation derived in lemma
 \ref{teo:mismasraices}. To start with, we work at the level of
 residual representations. Observe that the same relation proved in
 lemma \ref{teo:mismasraices} holds for the roots of the matrices in the image of the
 residual representations $\bar{\rho}_\lambda$ and
 $\bar{\sigma}_{1,\lambda}$:
 \begin{lema}
 \label{teo:primerencuentro}: Let $\lambda$ be a prime in $\Lambda_3$,
 then in the decomposition $\bar{\rho}_\lambda \cong \pi_{1,\lambda} \oplus
 \pi_{2,\lambda}$ we have $\pi_{1,\lambda} \cong
 \bar{\sigma}_{1,\lambda}$.
 \end{lema}
 Remark: Of course, we should write the above equality with $\pi_{i,\lambda}$
 for $i=1$ or $2$. But to fix notation, we will  always call $\pi_{1,\lambda}$
 the component of $\bar{\rho}_\lambda$ where the inertia group at
 $\ell$ acts with weights $0$ and $2k-3$ (as we did in section 3.1),
  this is a good way to
 distinguish the two components, and of course this is the only
 component that deserves being compared to $\bar{\sigma}_{1,\lambda}$.\\
Proof: Take $\lambda \in \Lambda_3$. Let $L$ be the Galois field corresponding to
$\bar{\rho}_\lambda$, thus $\Gal(L/ \mathbb{Q}) \cong M_\lambda$,
 and $B$ the one corresponding to
$\bar{\sigma}_{1,\lambda}$, thus if $\F_\lambda'$ is the residue field of
 $\mathcal{O}'_\lambda$ and $U_\lambda = \{ A \in \GL(2, \F_\lambda') :
 \det(A) \in \F_\ell^{2k-3}\}$,
 $\Gal(B/ \mathbb{Q}) \cong U_\lambda$. We want to prove that $B \subseteq L$.
Let $M = L \cap B$, and consider an element $z \in \Gal(B/M)$. Let $ \check{z}$
be a preimage of $z$ in $\Gal(\bar{\Q} / M)$, that we can choose such
that $\bar{\rho}_\lambda ( \check{z} ) = \mathbf{1}_4$ (because it is trivial on
$M = L \cap B$). Then, the residual version of lemma
 \ref{teo:mismasraices} implies that $1$ is a double root of the characteristic
 polynomial of $\bar{\sigma}_{1,\lambda}(\check{z})$. This implies that the group $\Gal(B/M)$
 is unipotent, but this group is a normal subgroup of $\Gal(B/\mathbb{Q}) \cong
  U_\lambda$, and $U_\lambda$ has no non-trivial unipotent normal
  subgroup, thus $B=M$, i.e., $B  \subseteq L$.\\
  Then, we have a projection: $\phi: \Gal(L/\mathbb{Q}) \rightarrow \Gal(B/
  \mathbb{Q})$, that is to say, $\phi$ sends $M_\lambda$ onto
  $U_\lambda$ and thus $\bar{\sigma}_{1,\lambda}$ is a quotient of $\bar{\rho}_{\lambda}$.\\
  Since $\bar{\rho}_{\lambda} \cong \pi_{1,\lambda} \oplus \pi_{2,\lambda}$ we conclude that
   $\bar{\sigma}_{1,\lambda} \cong \pi_{i, \lambda}$ with $i=1$ or $2$, and using the information on the Hodge-Tate decompositions
    we see that $i=1$.\\

\subsection{Proof of Theorem 2.1}
 We start by observing that part of the proof of lemma \ref{teo:primerencuentro}
  can be translated to the $\lambda$-adic setting. Take $\lambda \in
  \Lambda_3$, and call $L'$ the (infinite)
  Galois field corresponding to $\rho_\lambda$ and $B'$ the one
  corresponding to $\sigma_{1,\lambda}$. Recall that $\Gal(B' /\Q)$ is
  isomorphic to the subgroup $U'_\lambda$ of $\GL(2,
  \mathcal{O}'_\lambda)$ composed of matrices with determinant in
  $\mathbb{Z}_\ell^{2k-3}$, and therefore again we have a group with no
  non-trivial unipotent subgroups, thus we conclude from lemma \ref{teo:mismasraices}
   as in the proof of
  lemma \ref{teo:primerencuentro} that $B' \subseteq L'$ and that we
  have  a projection: $\phi': \Gal(L'/\mathbb{Q}) \rightarrow \Gal(B'/
  \mathbb{Q})$. We have $\phi' \circ \rho_\lambda = \sigma_{1,\lambda}$.
  Let us consider the normal subgroup $\Gal(L' /B')$ of $\Gal(L'/\Q)$,
  i.e., we are considering the  restriction $\rho_\lambda |_{\ker
  \phi'}$.
  The elements in this subgroup fix $B'$, thus by lemma \ref{teo:mismasraices}
   we see that the corresponding matrices in $\GSp(4,
   \mathcal{O}_\lambda)$ will have $1$ as a double root. \\
    On the other
   hand, we know that the residual representation $\bar{\rho}_\lambda \cong
   \pi_{1, \lambda} \oplus \pi_{2,\lambda} \cong \bar{\sigma}_{1,\lambda} \oplus
    \pi_{2, \lambda} $ has maximal image $M_\lambda$ (see lemmas
    \ref{teo:reducibleperogrande} and
    \ref{teo:primerencuentro}). Moreover,
     the representation $\sigma_{1,\lambda}$ being a
     ``deformation"
     of $\pi_{1,\lambda}$ which is disjoint from $\pi_{2,\lambda}$ (in the sense
      established during the proof of lemma \ref{teo:reducibleperogrande}, i.e.,
       up to the equality of determinants), we see that restricting to $\ker \phi'$
       will only shrink the image of $\pi_{2,\lambda}$ by making the
       determinant trivial, in other words:
        the residual representation
    $\overline
{\rho_\lambda |_{\ker \phi'}}$ has image $\SL(2, \F_{2,\lambda}) \oplus
    \mathbf{1}_2 \subseteq M_\lambda \quad \qquad (3.2)$. \\
    So, what do we know about $\rho_\lambda  |_{\ker \phi'} $? We have
    determined its residual image and we also know that all matrices in
    its image have $1$ as a double root: this last property extends to
    the Zariski closure of the image, and using the information we have
    together with the list of possibilities for this Zariski closure
    given in [T1], we see that the image of $\rho_\lambda  |_{\ker \phi'} $
     must be contained in $\SL(2, \mathcal{O}_\lambda) \oplus
     \mathbf{1}_2$. If we call $\mathcal{O}''_\lambda  \subseteq
     \mathcal{O}_\lambda$ the field generated by the traces of the
     image of
     $\rho_\lambda  |_{\ker \phi'} $, we can apply a lemma of Serre (cf. [S1])
      (and Carayol's lemma for the assertion about the field of definition, cf. [C])
       and
       conclude from (3.2) that the image of $\rho_\lambda  |_{\ker \phi'} $
      must be conjugated to  $\SL(2, \mathcal{O''}_\lambda) \oplus
     \mathbf{1}_2$.\\
     Remark: $\ker \phi'$ fixes $B'$ which is an infinite extension of $\Q$,
     but Serre's lemma can still be applied because $G_\Q$ is compact
     and the fixer of $B'$ is a closed subgroup.\\
     Thus, we conclude that $\Image (\rho_\lambda) \subseteq
      \GSp(4, \bar{E}_\lambda )$  contains a normal
     subgroup isomorphic to  $\SL(2, \mathcal{O''}_\lambda) \oplus
     \mathbf{1}_2$, and the quotient by this subgroup gives
     $U'_\lambda$. But it is easy to see that
     the normalizer of  $\SL(2, \mathcal{O''}_\lambda) \oplus
     \mathbf{1}_2$ in $\GSp(4, \bar{E}_\lambda )$ is
     contained in the reducible group $\GL(2,\bar{E}_\lambda )
      \oplus \GL(2, \bar{E}_\lambda )$. Thus, $\rho_\lambda$
      is reducible, for every $\lambda \in \Lambda_3$, and $\sigma_{1,\lambda}$
      is one of its two-dimensional irreducible components.

\section{From Theorem 2.1 to Theorem 1.1}

Recall that the results and proofs given in previous sections date from 2003. The result of ``existence of compatible families" that we proved in section 3.2 were extended in [D3], which is a ``sequel" to this paper, where it was applied to prove some cases of the Fontaine-Mazur conjecture.\\ Moreover, this result is key in the proof of Serre's conjecture given in [D4], [KW1], [K], [D5] and [KW2].\\
In the case of Hodge-Tate weights $\{ 0, 1 \}$ the results of potential modularity of Taylor do not imply that the representation is motivic, but as observed in [D3] we still can apply the techniques explained in section 3.2 and prove existence of families. In this case, the natural restriction becomes $\ell_0 > 2$.\\

As in the previous sections, we assume that we are not in the Saito-Kurokawa case (thus the reducible case must be a case of
 two $2$-dimensional irreducible components).\\

Thus, if $\ell_0$ is odd, $\ell_0 \nmid N$, we consider the irreducible component $\sigma_{2,\lambda_0}$ of $\rho_{\lambda_0}$ having
 Hodge-Tate weights $\{k-2 , k-1 \}$ and we apply existence of compatible families to $\sigma_{2,\lambda_0} \otimes \chi^{2-k}$ instead of $\sigma_{1,\lambda_0}$.\\
 The rest of the proof given in the previous sections extends word by word, except that $\sigma_{2,\lambda_0}$ and $\sigma_{1,\lambda_0}$ exchange roles. We conclude that theorem \ref{teo:unifow} is still true if we change the assumption $\ell_0 > 4k-5$ by $\ell_0 > 2$.\\

We have shown that for almost every prime $\lambda$ the representation $\rho_\lambda$ in our compatible family is reducible, and one of its $2$-dimensional irreducible components is $\sigma_{2,\lambda}$, a representation that lies in a strongly compatible family. It is obvious (from the definition of compatibility) that the second components, even if they are a priori defined only for almost every prime, will also form a compatible family, let us call them $\sigma_{1,\lambda}$. Moreover, from the formula:
$$ \rho_\lambda \cong \sigma_{1, \lambda} \oplus \sigma_{2, \lambda} $$
we see not only the compatibility of the $\sigma_{1,\lambda}$ but also that the representations $\sigma_{1, \lambda}$ are crystalline if $\ell$ is not in the ramification set $S$, of Hodge-Tate weights $\{0, 2k-3 \}$, with ramification set contained in $S$ and semistable or unramified locally at primes of $S$. \\
At this point, we apply an argument based on Serre's conjecture, which is now a theorem (cf. [D5] and [KW2]). Serre proved in [S2] using his conjecture that  compatible families of Galois representations as   $\{ \sigma_{2,\lambda}\otimes \chi^{2-k} \}$
 or $\{ \sigma_{1,\lambda} \}$ are modular. The fact that $\{ \sigma_{1,\lambda} \}$ is a priori only defined for almost every $\lambda$ is irrelevant for the argument of Serre: he only needs residual modularity in infinitely many characteristics, not in ALL characteristics (it is a typical patching argument). The essential condition is that the family is compatible, with constant Hodge-Tate weights and uniformly bounded conductor.\\
 We conclude that the families $\{ \sigma_{2,\lambda}\otimes \chi^{2-k} \}$ and $\{ \sigma_{1,\lambda} \}$ correspond to representations attached to classical modular forms, of weight $2$ and $2k-2$, respectively, and obviously this implies in particular that the family $\{ \sigma_{1,\lambda} \}$ is also defined for EVERY prime $\lambda$. This concludes the proof of theorem 
 \ref{teo:unifo}.\\
 
 Remark 1: Observe that in the particular case of a level $1$ Siegel cuspform, we conclude irreducibility for $p>2$ if it is not of Saito-Kurokawa type. The reason is that one of the $2$-dimensional components  would give rise to a level $1$ weight $2$ classical modular form.\\
 
 Remark 2: What we have shown is that our result of uniformity of reducibility, combined with Serre's conjecture (now a theorem), implies the truth of the Endoscopy Conjecture in the semistable case.\\
 
 Remark 3: Our result of ``uniformity of reducibility" can be extended to the non-semistable case, with the same arguments. In fact, this have been done recently by Skinner and Urban (cf. [SU], section 3). The argument of Serre explained above also applies here, so also the Endoscopy Conjecture follows in this case.

\section{Final Remarks}
 The Galois representations attached to a Siegel cusp form $f$ of level greater
 than one are
known to verify the semistability condition when the ramified
local components of (the automorphic representation corresponding
to) $f$ are of certain particular types (for example, a Steinberg
representation), as follows from  recent works of
Genestier-Tilouine and Genestier (cf. [GT]). Thus, the results in this
article apply to these cases.
 We thank J. Tilouine for pointing out this fact to us. \\


\section{Bibliography}

[C] Carayol, H., {\it  Formes modulaires et repr{\'e}sentations galoisiennes {\`a}
 valeurs dans un anneau local complet}, Contemp. Math. {\bf 165}, Amer.
  Math. Soc. (1994) 213-237

[D1] Dieulefait, L., {\it On the images of the Galois representations
 attached to genus $2$
Siegel modular forms}, J. Reine Angew. Math. {\bf 553} (2002) 183-200

[D2] Dieulefait, L.,  {\it Irreducibility of Galois actions on
level $1$ Siegel cusp forms}, Modular Curves and Abelian Varieties, Cremona et al. (eds.), Birkhauser,
 Progress in Math. {\bf 224} (2004) 75-80

[D3] Dieulefait, L., {\it Existence of families of Galois
representations and new cases of the Fontaine-Mazur conjecture}, J. Reine Angew. Math. {\bf 577} (2004) 147-151

[D4] Dieulefait, L., {\it The level 1 weight 2 case of Serre's conjecture}, to appear

[D5] Dieulefait. L., {\it Remarks on Serre's modularity conjecture}, preprint available at: www.arxiv.org

[GT] Genestier, A., Tilouine, J., {\it Syst\`{e}mes de Taylor-Wiles pour $\GSp_4$}, in
 Formes Automorphes (II), Asterisque {\bf 302} (2005) 177-290

[K] Khare, C., {\it Serre's modularity conjecture: the level one case}, Duke
Math. J. {\bf 134} (2006) 557-589

[KW1] Khare, C., Wintenberger, J-P., {\it On Serre's reciprocity conjecture for $2$-dimensional mod $p$ representations of the Galois group of $\Q$}, preprint, (2004); available at: www.arxiv.org

[KW2] Khare, C., Wintenberger, J-P., {\it Serre's modularity conjecture (I, II)}, preprint, (2006);
 available at: \\
 http://www.math.utah.edu/~shekhar/papers.html

[R] Ribet, K., {\it Images of semistable Galois
representations}, Pacific J.  Math. {\bf 181} (1997)

[S1] Serre, J-P.,  {\it Abelian $\ell$-adic representations and elliptic
curves}, (1968) Benjamin

[S2] Serre, J-P., {\it Sur les repr{\'e}sentations modulaires de degr{\'e}
$2$ de $\Gal(\bar{\mathbb{Q}} / \mathbb{Q})$}, Duke Math. J. {\bf 54}
(1987) 179-230

[SU] Skinner, C., Urban, E. {\it Sur les def\'{o}rmations p-adiques de certaines repr\'{e}sentations automorphes}, J. Inst. Math. Jussieu {\bf 5} (2006) 629-698 

[T1] Taylor, R., {\it On the $\ell$-adic cohomology of Siegel
threefolds}, Invent. Math. {\bf 114} (1993) 289-310

[T2] Taylor, R., {\it Remarks on a conjecture of Fontaine and Mazur},
J. Inst. Math. Jussieu {\bf 1} (2002)

[T3] Taylor, R., {\it On the meromorphic continuation of degree two
 L-functions}, preprint, (2001);
 available at http://abel.math.harvard.edu/$\sim$rtaylor/

[T4] Taylor, R., {\it Galois Representations}, Annales Fac. Sc. Toulouse {\bf 13} (2004) 73-119 

[W1] Weissauer, {\it The Ramanujan Conjecture for genus two
Siegel modular forms (An application of the Trace Formula)}, preprint (1994)

[W2] Weissauer, {\it Four dimensional Galois representations}, in Formes Automorphes (II), Asterisque {\bf 302} (2005) 67-150

\end{document}